# Walking in the OEIS:
# From Motzkin numbers to Fibonacci numbers.
# The «shadows» of Motzkin numbers

Gennady Eremin

ergenns@gmail.com

August 24, 2021

**Abstract.** In this paper, we consider nine OEIS sequences, the analysis of which allows us to find a connection between Motzkin numbers and Fibonacci numbers. In each Motzkin number, we distinguish an even component and an odd component, the difference of these two components is called the "shadow" of the Motzkin number. Reverse of the Motzkin shadows give us a sequence from the family of Fibonacci numbers.

*Keywords*: On-Line Encyclopedia of Integer Sequences (OEIS), Motzkin paths, Motzkin numbers, Fibonacci numbers, generating function, series reversion.

In this paper, we will work with the On-line Encyclopedia of Integer Sequences (OEIS) [Slo21]. On the Internet, all sequences are available to the reader; there is a separate website for each sequence. Often in formulas borrowed from various sites, e.g., generating functions in closed form (explicit formulas), we will retain the notation adopted in OEIS. For example, * is a multiplication, ^ is a power, sqrt is the square root, etc.

We consider nine OEIS sequences: A000045, A001006, A007440, A039834, A100223, A107587, A214649, A343386, and A343773, the analysis of which gives us the opportunity to establish a connection between Motzkin numbers and Fibonacci numbers. It all started by identifying an even component and an odd component in each Motzkin number.

## 1    Introduction

Even and odd elements are often distinguished in integer sequences. In OEIS, the most well-known case is the division of a sequence of natural numbers A000027 into two subsequences: a sequence of even elements A005843 (with zero) and a sequence of odd elements A005408. In addition, there is a sequence A299174 – even positive numbers (without zero).

Another example, the well-known Catalan numbers A000108 are divided into even and odd indices. As a result, we have a sequence of Catalan numbers with even indexes A048990 and a sequence with odd indexes A024492. Additionally, we note that in the sequence A048990, all terms are even, so all odd Catalan numbers fall into the sequence A024492. Also, Emeric Deutsch and Bruce Sagan proved that in A000108 the Catalan numbers with indices $n = 2^k - 1$, $k \geq 0$, are odd [DS04].

Finally, we note the famous Indian mathematician H. L. Manocha, who obtained generating functions for certain polynomials by splitting a given series into two parts with even and odd elements [Man74]. In OEIS, there are hundreds of thousands of sequences, and the reader can give additional examples where the elements of a power series are separated by parity.

In this paper, we are working with Motzkin numbers, and we will not distribute these numbers over different sequences. We divide each Motzkin number into two parts – the so-called even component and the odd one, and these components form two new sequences.

**Sequence A001006**. We show the beginning of the sequence of Motzkin numbers $m_i$, $i \geq 0$:

(1.1)     1, 1, 2, 4, 9, 21, 51, 127, 323, 835, 2188, 5798, 15511, 41835, ...

On the website A001006, the reader will find a convenient recurrence relation

(1.2)     $(n+2)\, m_n = (2n+1)\, m_{n-1} + 3(n-1)\, m_{n-2}, \ n \geq 2.$

*The generating function* for an arbitrary sequence of integers $g_i$, $i \geq 0$, is a formal power series

$$G(x) = g_0 + g_1 x + g_2 x^2 + \ldots + g_n x^n + \ldots$$

Accordance to (1,1), the generating function for Motzkin numbers is

(1.3)     $M(x) = 1 + x + 2x^2 + 4x^3 + 9x^4 + 21x^5 + 51x^6 + 127x^7 + \ldots$

Motzkin numbers count Motzkin words, which can be defined in various ways. In (1.1), the $n$th element is equal to the number of words of length $n$. The elements are indexed from zero, i.e., a Motzkin word of length 0 (an empty word) is allowed.

Often work with correct bracketed sets with zeros. In our case, we are dealing with Motzkin paths on a square lattice. Each path includes steps of only three types: diagonal up-step U = (1, 1), diagonal down-step D = (1, –1) and horizontal step H = (1, 0). A path of length $n$ starts at the origin, ends at the point $(n, 0)$ and never goes below the $x$-axis. Motzkin paths are words on the alphabet {U, D, H}. In each Motzkin path, the number of steps U and D is always the same. A Motzkin path can include only horizontal steps, and it can also be without horizontal steps. A Motzkin path with no H steps is called a Dyck path.

In OEIS sequences, one can find an equation for generating functions. You can get the equation for $M(x)$ using the *inference rules* of the Motzkin language. A Motzkin word is: (1) an *empty word* λ, or (2) a word of the form H$a$, where $a$ is a Motzkin word, or (3) a word of the form U$a$D$b$, where $a$ and $b$ are Motzkin words.

Let $\mathcal{M}$ be the set of Motzkin words, then the listed three rules correspond to the structural equation

(1.4) $\quad\quad\quad\quad\quad \mathcal{M} = \lambda + H\mathcal{M} + U\mathcal{M}D\mathcal{M}$,

where plus denotes the union of disjoint sets. Replacing the set $\mathcal{M}$ with the generating function $M(x)$, we get the *functional equation* [DS77]

(1.5) $\quad\quad\quad\quad\quad M(x) = 1 + xM(x) + x^2 M^2(x)$.

In (1.4), we replaced each character U, D and H with $x$, and the only empty word $\lambda$ changes to 1. The solution of the quadratic equation (1.5) gives us the following expression in closed form for the generating function (recall that in such explicit formulas we will keep the notation adopted in OEIS):

(1.6) $\quad\quad\quad$ M (x) = ( 1 – x – sqrt (1 – 2*x – 3*x^2) ) / (2* x^2).

## 2   Even and odd Motzkin paths

By the number of U steps (or D steps), we determine whether the path is even or odd. Let's formulate the corresponding definition.

**Definition 1**. *A Motzkin path is called even if it has an even number of U steps. Otherwise, the Motzkin path is odd.*

Like many mathematicians, we consider zero to be an even number, so we consider a path in which there are only horizontal steps to be even. For example, of the nine Motzkin paths of length 4, three paths are even – these are HHHH, UUDD, and UDUD. Accordingly, the number of odd paths of length 4 is 9 – 3 = 6. Let's list them: HHUD, HUHD, HUDH, UHHD, UHDH, and UDHH. Obviously, we have one even Motzkin path of length from 0 to 3 inclusive – this is an empty word $\lambda$ and the paths H, HH, HHH. Let's add that there are no odd Motzkin paths of length 0 or 1, and there is only one UD path of length 2.

In OEIS, there are two corresponding sequences A107587 and A343386. The sequence A107587 was created in 2005, but recently significantly corrected and supplemented by the author together with Sergey Kirgizov (LIB, Univ. Bourgogne Franche-Comté, France). The second sequence A343386 was placed in OEIS this year by the same authors.

**Sequence A107587** contains numbers that count the number of even Motzkin paths; here is the beginning of this sequence

(2.1) $\quad\quad$ 1, 1, 1, 1, 3, 11, 31, 71, 155, 379, 1051, 2971, 8053, 21165, ...

The terms of the sequence A107587 are indexed from zero (like the Motzkin numbers in A001006). The generating function for A107587 has the form

(2.2) $\quad A(x) = 1 + x + x^2 + x^3 + 3x^4 + 11x^5 + 31x^6 + 71x^7 + \ldots$

**Proposition 2.** *All terms of the sequence A107587 are odd.*

*Proof.* On the website A107587, the reader will find the author's formula for the terms of the sequence; the formula uses binomial coefficients and Catalan numbers:

$$\sum_{k \geq 0} \binom{n}{4k} \text{Cat}(2k), \ n \geq 0.$$

In such a sum, for k = 0, the initial number is 1, the remaining summands are even due to the parity of the Catalan numbers with even indices 2, 4, 6,... Thus, the value of the sum is always an odd number. □

**Sequence A343386** contains numbers that count the number of odd Motzkin paths; here is the beginning of this sequence:

(2.3) $\quad$ 0, 0, 1, 3, 6, 10, 20, 56, 168, 456, 1137, 2827, 7458, 20670, ...

The elements of the sequence (2.3) are also indexed from zero. The generating function for A343386 can be written as

(2.4) $\quad B(x) = x^2 + 3x^3 + 6x^4 + 10x^5 + 20x^6 + 56x^7 + 168x^8 + \ldots$

The following functional equation is obvious:

$$A(x) + B(x) = M(x).$$

The similar amount is valid for terms

$$A107587(n) + A343386(n) = A001006(n), \ n \geq 0.$$

**Generating functions $A(x)$ and $B(x)$.** Let's derive the equations for A(x) and B(x). We formulate the inference rules for even and odd Motzkin paths. An even Motzkin path is:

  – an empty path λ, or
  – a path of the form H*a*, where *a* is an even Motzkin path, or
  – a path of the form U*a*D*b* or U*b*D*a*, where *a* is an even Motzkin path, and *b* is an odd Motzkin path.

An odd Motzkin path is:

- a path of the form H*b*, where *b* is an odd Motzkin path, or
- a path of the form U*a*D*a'*, where *a* and *a'* are even Motzkin paths, or
- a path of the form U*b*D*b'*, where *b* and *b'* are odd Motzkin paths.

Let's denote by $\mathcal{A}$ и $\mathcal{B}$, respectively, the sets of even and odd Motzkin paths. Then the following two structural equations correspond to the listed inference rules:

(2.5) $$\mathcal{A} = \lambda + H\mathcal{A} + U\mathcal{A}D\mathcal{B} + U\mathcal{B}D\mathcal{A},$$

(2.6) $$\mathcal{B} = H\mathcal{B} + U\mathcal{A}D\mathcal{A} + U\mathcal{B}D\mathcal{B}.$$

Still, plus denotes the union of disjoint sets. After replacing the set $\mathcal{A}$ with the generating function $A(x)$, and the set $\mathcal{B}$ with the generating function $B(x)$, we get the functional equalities

(2.7) $$A(x) = 1 + xA(x) + 2x^2A(x)B(x),$$

(2.8) $$B(x) = xB(x) + x^2A^2(x) + x^2B^2(x).$$

We have again replaced each character U, D and H with $x$, and the empty word $\lambda$ is changed by 1. It is easy to see, the addition of equations (2.7) and (2.8) gives us the functional equation of the generating function for Motzkin numbers, i.e., equality (1.5).

If we replace $B(x) = M(x) - A(x)$ in formula (2.7), and then use formula (1.6), we get the quadratic equation

(2.9) $$2*x^2*A(x)^2 + \text{sqrt}(1-2*x-3*x^2)*A(x) - 1 = 0.$$

The solution of (2.9) gives a closed-form formula for $A(x)$ (see A107587):

(2.10) $$A(x) = (\text{sqrt}(1 - 2*x + 5*x^2) - \text{sqrt}(1 - 2*x - 3*x^2)) / (4*x^2).$$

It is not difficult to obtain a closed-form expression for $B(x)$, it is enough to subtract equality (2.10) from equality (1.6):

(2.11) $$B(x) = (2 - 2*x - \text{sqrt}(1-2*x-3*x^2) - \text{sqrt}(1-2*x+5*x^2)) / (4*x^2).$$

We have divided the *n*th Motzkin number into even and odd components, which count, respectively, the number of even and odd Motzkin paths of length *n*. Interesting is not only the sum of these two components (which is equal to the original Motzkin number), but also their difference. Next, we will show that a sequence composed of such difference numbers brings us closer to the Fibonacci numbers.

# 3  Shadows of Motzkin numbers

The author and Sergey Kirgizov designed another OEIS sequence A343773 – the excess of the number of even Motzkin paths over the number of odd Motzkin paths. That is,

$$A107587(n) - A343386(n) = A343773(n), \quad n \geq 0.$$

Here is the beginning of the sequence A343773:

(3.1)  1, 1, 0, –2, –3, 1, 11, 15, –13, –77, –86, 144, 595, 495, ...

As we can see, in (3.1) there are negative numbers among the terms s(n), n ≥ 0. On the website A343773, the reader will find the following recurrence relation:

(3.2)  $s(n) = ((2*n+1)*s(n-1) - 5*(n-1)*s(n-2)) / (n+2), \quad n > 1.$

In accordance with (3.1), the generating function, denoted by $S(x)$, has the form

(3.3)  $S(x) = 1 + x - 2x^3 - 3x^4 + x^5 + 11x^6 + 15x^7 - 13x^8 - 77x^9 - \ldots$

In order to obtain a functional equation for $S(x)$, subtract equality (2.8) from (2.7).

$$\begin{aligned}S(x) &= A(x) - B(x) \\ &= (1 + xA(x) + 2x^2 A(x) B(x)) - (xB(x) + x^2 A^2(x) + x^2 B^2(x)) \\ &= 1 + xS(x) - x^2 S^2(x).\end{aligned}$$

As a result, we get

(3.4)  $S(x) = 1 + xS(x) - x^2 S^2(x).$

The functional equation for $S(x)$ practically coincides with the corresponding equation of the generating function for Motzkin numbers (1.5), the difference is only in the last sign. We have called the terms of the sequence A343773 the *shadows* of Motzkin numbers.

An explicit formula of the generating function for the shadows of Motzkin numbers can be obtained by subtracting the equality (2.11) from (2.10):

(3.5)  $S(x) = (-1 + x - \sqrt{1 - 2x + 5x^2}) / (2x^2).$

In OEIS, there are duplicates of A343773. These are two sequences A100223 and A214649, differing only in their initial terms. There is also the sequence A007440, in which some of the terms have other signs. Let's consider each of these sequences separately.

**Sequence A100223** contains all the terms from A343773 plus two additional initial terms. Let's show the beginning of A100223:

(3.6)    1, 0, 1, 1, 0, –2, –3, 1, 11, 15, –13, –77, –86, 144, 595, 495, ...

The elements in (3.6) are indexed from zero, and we can write an obvious formula connecting the terms of both sequences:

$$A343773(n) = A100223(n+2), \ n \geq 0.$$

It is not difficult to obtain a functional equation of the generating function for the sequence A100223:

$$1 + x^2 + x^3 - 2x^5 - 3x^6 + x^7 + 11x^8 + 15x^9 - 13x^{10} - 77x^{11} - \ldots = 1 + x^2 S(x).$$

**Sequence A214649** differs from A100223 only by the second term:

(3.7)    1, –1, 1, 1, 0, –2, –3, 1, 11, 15, –13, –77, –86, 144, 595, 495, ...

Unlike (3.6), the elements in (3.7) are indexed starting from –1. The formula connecting the terms of the sequences A343773 and A214649 is obvious:

$$A343773(n) = A214649(n+1), \ n \geq 0.$$

Let's write a functional equation of the generating function for A214649:

$$1/x - 1 + x + x^2 - 2x^4 - 3x^5 + x^6 + 11x^7 + 15x^8 - 13x^9 - \ldots = 1/x - 1 + xS(x).$$

We are most interested in the third sequence A007440, which is called "Reverse of the generating function for Fibonacci numbers" and which differs from A343773 only in the signs of some terms.

**Sequence A007440,** unlike A100223 and A214649, does not have additional terms. If you do not take into account the signs of numbers, the sequences A343773 and A007440 coincide. Here is the beginning of A007440:

(3.8)    1, –1, 0, 2, –3, –1, 11, –15, –13, 77, –86, –144, 595, – 495, ...

Additionally, we note that the elements in A007440 are indexed from 1.

Let's compare the lists (3.1) and (3.8). It is easy to see, every number in (3.8) with an odd index of 1, 3, 5, etc. retains the sign in (3.1). While the elements with even indexes change sign to the opposite. The terms in A343773 are indexed from 0, and we can write the following formula:

(3.9)    $A343773(n) = (-1)^n \times A007440(n+1), \ n \geq 0.$

In equality (3.9), we used a well-known procedure called *inversion*. The inversion is symmetric and often occurs in OEIS. In our case, the sequence A343773 is the inversion of A007440, and vice versa, A007440 is the inversion of A343773, i.e.,

$$\text{A007440}(n+1) = (-1)^n \times \text{A343773}(n), \quad n \geq 0.$$

Let's write the equation of the generating function for A007440:

$$x - x^2 + 2x^4 - 3x^5 - x^6 + 11x^7 - 15x^8 - 13x^9 + 77x^{10} - 86x^{11} - \ldots = xS(-x).$$

As we can see, in the functional equation, the inversion procedure is reduced to replacing the sign of a part of the terms in the generating function. Most often, inversion is used precisely to change signs. Note that the inversion procedure is universal; there are no restrictions on the processed series.

Another well-known procedure is called the *reverse* of the generating function. It is the reverse procedure that is mentioned in the title of A007440 in connection with Fibonacci numbers. Usually, the reverse radically changes the given sequence. For example, it is enough to compare the sequences A007440 and A000045.

**From Motzkin numbers to Fibonacci numbers.** Figure 1 shows the steps that we have already taken on the way from Motzkin numbers to Fibonacci numbers. First, we divided an arbitrary Motzkin number (the term of A001006) into an even component (the term of A107587) and an odd component (the term of A343386). Then the difference of these components is calculated, and we get the shadow of the Motzkin number (the term of A343773). The sequence A343773 is duplicated in A100223 and in A214649, which are also shown in the figure. We can say that three sequences A100223, A214649 and A343773 form a cluster of duplicates.

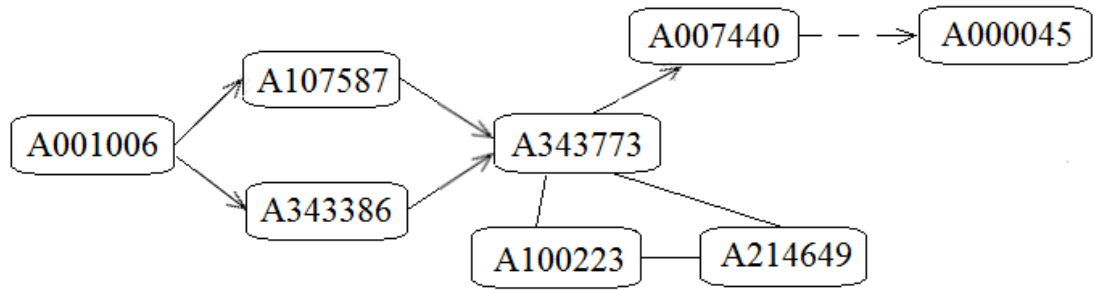

Figure 1: Step-by-step technological chain.

At the next step, the term of A343773 is converted into an element of A007440 using the inversion procedure. And finally, using the reverse procedure, we get the Fibonacci number (the term of A000045). However, the last transformation has not yet been verified by us, so a broken line leads from the sequence A007440 to A000045.

Below we will consider the reverse procedure, construct a reverse power series for the shadows of Motzkin numbers, which, as it turns out, belongs to the family of Fibonacci numbers.

## 4 Series reversion

There is no generally accepted definition of the reverse power series for a given integer sequence in the literature. Many authors, describing the reverse procedure, are guided by the classical work edited by M. Abramowitz and I. A. Stegun (see [AS72], p. 16), where the corresponding algorithm is described in sufficient detail. Additionally, we note the encyclopedic network source Wolfram Web Resource [Wei21].

The algorithm of the series reversion is reduced to calculating the coefficients of the inverse generating function, taking into account the coefficients of a given generating function. So, we take the original power series without a free term (the constant term is zero, the elements of the series are indexed from 1):

(4.1) $$y = ax + bx^2 + cx^3 + dx^4 + ex^5 + fx^6 + \ldots$$

Starting from the coefficients of the sequence (4.1), we need to obtain the coefficients of the reverse series

(4.2) $$x = Ay + By^2 + Cy^3 + Dy^4 + Ey^5 + Fy^6 + \ldots$$

Let's substitute the sum (4.2) instead of $x$ in each term of (4.1), and then we get the following equality:

(4.3) $$y = (aA)y + (aB + bA^2)y^2 + (aC + 2bAB + cA^3)y^3 + \ldots$$

It remains to equate in (4.3) the coefficients with the same degrees of $y$ on both sides of the equality:

(4.3a)    $1 = aA$, or $A = a^{-1}$;
(4.3b)    $0 = aB + bA^2 = aB + ba^{-2}$, then $B = -ba^{-3}$;
(4.3c)    $0 = aC + 2bAB + cA^3 = aC + 2ba^{-1}(-ba^{-3}) + ca^{-3}$, |
         then $C = 2b^2a^{-5} - ca^{-4}$.

Here is another coefficient:

(4.3d)    $D = 5bca^{-6} - da^{-5} - 5b^3a^{-7}$.

It is not difficult to see, the calculation of further coefficients is very tedious. The procedure for reversing a power series is implemented in the Wolfram Language. A convenient software service allows the reader to perform fairly complex mathematical calculations online [WM21].

For some generating function $G(x)$, we denote by $G(x)^{rev}$ the reverse generating function. The reverse procedure, like the inversion procedure, is symmetric, so

$$(G(x)^{rev})^{rev} = G(x).$$

**Reverse of the shadows of Motzkin numbers.** We need to get the reverse series for A343773. Recall that the reverse procedure can be performed for a sequence that does not have a free term. In addition, from the formulas (4.3 a) – (4.3 d) it can be see that in the given series (4.1) the coefficient at $x$ is not equal to zero. Therefore, the indexing of the terms of the original series should start with 1. Thus, the sequence of shadows of Motzkin numbers (3.1) and, accordingly, the generating function (3.3) must be shifted to the right by one element. As a result, we will work with the generating function

(4.4) $\quad xS(x) = x + x^2 + 0x^3 - 2x^4 - 3x^5 + x^6 + 11x^7 + 15x^8 - 13x^9 - 77x^{10} - \ldots$

Let's calculate the first coefficients of the reverse series for (4.4) using the formulas (4.3 a) – (4.3 d); in our case $a = b = 1$, $c = 0$, $d = -2$.

$A = 1^{-1} = 1;$
$B = -1 \times 1^{-3} = -1;$
$C = 2 \times 1^2 \times 1^{-5} - 0 \times 1^{-4} = 2;$
$D = 5 \times 1 \times 0 \times 1^{-6} - (-2) \times 1^{-5} - 5 \times 1^3 \times 1^{-7} = 0 + 2 - 5 = -3.$

We show the beginning of the resulting reverse series:

(4.5)  1, –1, 2, –3, 5, –8, 13, –21, 34, –55, 89, –144, 233, –377, 610, –987, …

In (4.5), the terms $t(n)$ are indexed from 1, and the recurrent formula resembles the formula of Fibonacci numbers:

(4.6) $\quad\quad\quad\quad\quad\quad t(n) = -t(n-1) + t(n-2), \ n > 1.$

The generating function for the sequence (4.5) can be written as:

(4.7) $\quad (xS(x))^{rev} = x - x^2 + 2x^3 - 3x^4 + 5x^5 - 8x^6 + 13x^7 - 21x^8 + 34x^9 - \ldots$

Recall that the reverse is a symmetric procedure. The reader can calculate the initial elements of the reverse series for (4.7) using the formulas (4.3 a) – (4.3 d); in this case, $a = 1$, $b = -1$, $c = 2$, $d = -3$. As a result, the initial terms of A343773 will be obtained.

The sequence (4.5) and the corresponding generating function (4.7) practically coincide with the signed Fibonacci numbers A039834:

(4.8)  1, 1, 0, 1, –1, 2, –3, 5, –8, 13, –21, 34, –55, 89, –144, 233, –377, 610, …

In (4.8), the terms are indexed with -2, so the generating function for A039834 can be written as:

$$x^{-2} + x^{-1} + (x\,S(x))^{rev}.$$

Finally, to get to the Fibonacci numbers, we need to take the last step. There are no negative numbers in the Fibonacci sequence, and in A039834, every even term starting from zero is negative (numbers with odd indices are positive). Therefore,

(4.9) $$A000045(n) = -(-1)^n \times A039834(n), \ n \geq 0.$$

Obviously, equality (4.9) is equivalent to the inversion procedure.

Recall that we can get the Fibonacci numbers in a different way if we perform an inversion of the shadows of the Motzkin numbers, and then reverse the resulting sequence A007440.

## 5 The technological chain from Motzkin numbers to Fibonacci numbers

The final technological chain of the transition from the Motzkin numbers to the Fibonacci numbers is shown below in Figure 2; the transition is performed in four steps. In the first step, we divide the Motzkin numbers into an even component and an odd one, splitting the sequence A001006 into two sequences A107587 and A343386. In the second step, the difference of the new sequences gives us the shadows of the Motzkin numbers, A343773.

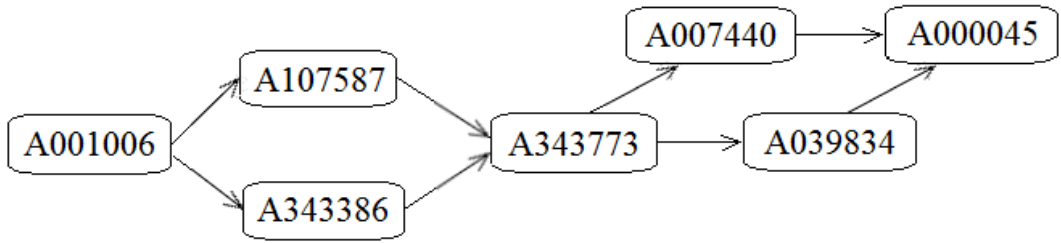

Figure 2: From Motzkin numbers to Fibonacci numbers.

Central place in the process chain is occupied by shadows of Motzkin numbers, starting from which we can obtain Fibonacci numbers in two ways: (1) we invert the shadows of Motzkin numbers, then we reverse the obtained sequence A007440, (2) we reverse the shadows of Motzkin numbers, then we invert the obtained sequence A039834.

Due to the symmetry of the inversion and reverse procedures, the last two steps of the technological chain can be performed both in the forward and back direction. This means that we can perform the transition back from the Fibonacci numbers to the shadows of the Motzkin numbers, and also in two ways: (1) invert the Fibonacci numbers, then reverse the resulting sequence A039834, (2) reverse the Fibonacci number, then invert the resulting sequence A007440.

**Acknowledgements.** I would like to thank Sergey Kirgizov (LIB, Univ. Bourgogne Franche-Comté, France) for his help in calculating the generating functions.

Gzhel State University, Moscow, 140155, Russia
http://www.art-gzhel.ru/